# PORTFOLIO SELECTION BY THE MEANS OF CUCKOO OPTIMIZATION ALGORITHM

Elham Shadkam, Reza Delavari, FarzadMemariani, MortezaPoursaleh

Department of Industrial Engineering, Faculty of Eng., Khayyam University, Mashhad, Iran

### ABSTRACT

*Portfolio selection is one of the most important and vital decisions that a real or legal person, who invests in stock market should make. The main purpose of this paper is the determination of the optimal portfolio with regard to stock returns of companies, which are active in Tehran's stock market. For achieving this purpose, annual statistics of companies' stocks since Farvardin 1387 until Esfand 1392 have been used. For analyzing statistics, information of companies' stocks, the Cuckoo Optimization Algorithm (COA) and Knapsack Problem have been used with the aim of increasing the total return, in order to form a financial portfolio. At last, results merits of choosing the optimal portfolio using the COA rather than Genetic Algorithm are given.*

## 1.INTRODUCTION

The Stock Market of Tehran is a dual purpose organization, which is a centre to collect savings and money of the private sector in order to finance the long-period investment projects and is a safe place to invest their surplus funds in different companies [1].

Nowadays, stock market acts a very important role in economic development using the tools like pricing, reducing the risk, resource mobilization and optimal allocation of capital [2].Selecting the optimal portfolio is the most important issue. There are different methods for selecting the optimal portfolio based on the results of the research in this field [3]. During the time, the defects of each method identified and new algorithm has been replaced.

One of the greatest theories in determining the optimal portfolio in past decades is the "Modern Portfolio Theory" which is presented by Harry Markowitz and William Sharpe. The Modern Portfolio Theory is a holistic approach to the stock market. Unlike the technical ones, this method discuss about the whole stock in market. In other words, this theory has a macro perspective versus a micro view. So the portfolio and the optimal combination of stocks are emphasized. Although in making a portfolio, the relationship between risk and return of stocks is important. One of the most important criteria of decision making in stock market is the return of stocks. Return of stock itself includes information and investors can use them in their financial analysis [4].

International Journal on Computational Sciences & Applications (IJCSA) Vol.5, No.3, June 2015

Also, in different studies about determining the optimal portfolio, "Risk" is known as one of the main criteria for determining the portfolio, which is clearly observable in earlier theories of Markowitz and other classic economists. Modern portfolio theory is useful in making the portfolios that have lower risk with respect to most return [5].

In the early 1950s Markowitz made the portfolio to quantity by defining the expected return as the average of returns and the risk as its variance. In his expanded model, investors can decrease the risk of their portfolio for certain return or increase the expected return of their portfolio for certain risk level.

Meta-heuristic algorithms have been used in portfolio selection problem including Genetic algorithm [6], [7], [8], [9]; particle swarm [10]; Ant colony algorithm [11], [12]; simulated annealing algorithm [13].

In this research some of the superior companies of stock market have been chosen and the information of their return is gathered. First, using the Knapsack model, portfolio selection problem changed into a mathematical model and then in order to select the optimal portfolio, Cuckoo optimization Algorithm (COA) is used for solving the Knapsack model. At last the Knapsack model is solved using the Genetic algorithm and the results compare with COA.

In the second section, COA algorithm has been introduced briefly. The third section introduces the data used in this article. The fourth section provides the results of implementing the COA in portfolio selection problem, proves the effectiveness of this algorithm, and draws some comparisons to a Genetic Algorithm implementation.

## 2. INTRODUCING THE CUCKOO OPTIMIZATION ALGORITHM

The Cuckoo optimization algorithm has been presented by Yang and Deb in 2009 [14]. This algorithm is inspired by the life of cuckoos in egg laying method and the levy flight instead of simple random isotropic step. Later, in 2011, the COA has been investigated with more details by Rajabioun [15].

Like other evolutionary algorithms, this one begins with an initial population of answers (cuckoos). These cuckoos have some eggs and they put these eggs in nests of other birds' and wait until the host bird maintains these eggs beside her eggs. In fact, this lazy bird forces other birds to survive her generation so nicely. Some of the eggs that have less similarity to the host bird's eggs will be recognized and destroyed.

In fact the cuckoos improve to imitate the target host bird's eggs and also the host birds learn how to recognize the fake eggs continuously.The number of survived eggs in each zone shows suitability of this zone, and the greater number of survived eggs, the more attention pays to that zone. In fact, this is the parameter that the COA wants to optimize. This algorithm applied to many problem and good results attained [16], [17]. The COA diagram is given Figure 1. In order to solve an optimization problem, the variable values of the problem should take shape out of an array called "habitat". In an optimization problem, the $N_{var}$ of a habitat will be an $1 \times N_{var}$ array that shows the current living location of cuckoos. This array describes as:

$$habitat = [X_1, X_2, ..., X_{N_{var}}]$$

The suitability (profit) in the current habitat obtains by evaluating the profit function ($f_p$) in the habitat. So:

$$Profit = f_b(habitat) = f_b(x_1, x_2, ..., x_{Nvar})$$

For starting an optimization algorithm, a habitat matrix in size of $N_{pop} * N_{var}$ is made, then for each habitat, a random number of eggs will be allocated. Considering the number of eggs that each cuckoo egg lays and also the distance between the cuckoos and the current optimized zone, the laying radius will be calculated and the cuckoos begin to egg lay in that zone.

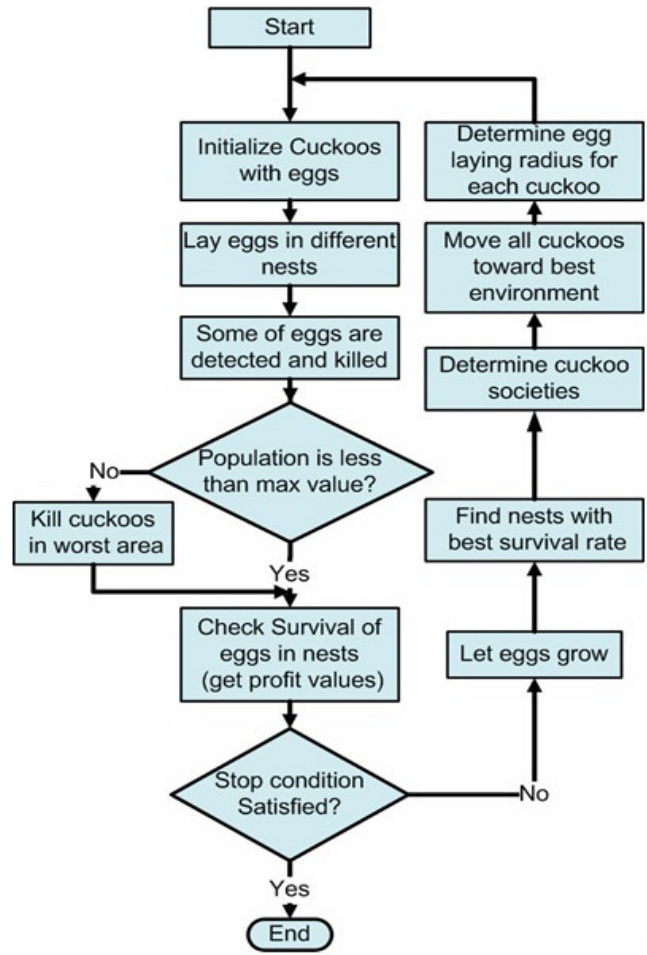

<Figure 1> Flowchart of COA [15]

The Egg Laying Radius (ELR) is given by

$$ELR = a \times \frac{Number\ of\ current\ cuckoo's\ eggs}{Total\ number\ of\ eggs} \times (Var_{hi} - Var_{low})$$

Then each cuckoo begins to egg laying in the nests that are in her ELR. So after each egg laying, p% of eggs (usually 10%) that is less profitable destroys. Other chicks grow up in the host nests.

## 2.1. THE CUCKOO'S MIGRATION

When the little cuckoos grow up and get mature, they live in their zone and groups for a while, but when the laying time comes they migrate to better zone that have a better possibility of survival.

After the cuckoos groups in general different living zones composed (justified region or search space of the problem), the group with the best location will be the migration target for other cuckoos

This is hard to find out each cuckoo belongs to which group when the grown cuckoos live all around the environment. For solving this problem, the cuckoos will cluster by the "K-means" method, which is a classic way to group (finding a K between 3 and 5 is usually acceptable).
As it is shown in Figure 2, when the cuckoos migrate to the target, they don't travel the direct way. They just travel a part of the way (λ %) also in that part there is deflection (φ).

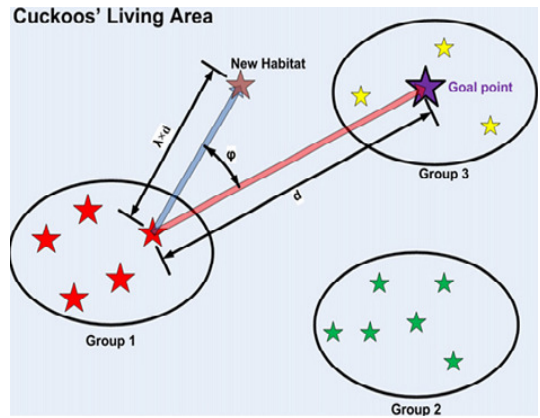

<Figure 2> The cuckoos' migration to the target [15]

These two parameters *(λ, φ)* helps cuckoos to explore a larger area. λ is a random number between 0 and 1 and φ is a number between $-\dfrac{\pi}{6}$ and $\dfrac{\pi}{6}$.

The migration formula is:

$$X_{\text{Next Habitat}} = X_{\text{current Habitat}} + F(X_{\text{Goal Point}} - X_{\text{current Habitat}})$$

## 2.2. THE CONVERGENCE OF ALGORITHM

After some iteration, all of the cuckoos migrate to an optimized point with the most similarity of eggs to the host ones and access to the richest food sources.

This location has the top profit and in this location we have the least number of the killed eggs. The convergence of more than 95% of all cuckoos to a single point put the COA at its end.

## 3. MODELING THE PORTFOLIO SELECTION PROBLEM USING KNAPSACK

The statistical population used in this paper includes all of the accepted and active companies in Tehran's stock market. From these companies some of them have been chosen to study on because they had the following qualifications:

1. Companies listed in Tehran's stock market since Farvardin's 1387 until Esfand's 1392 with uninterrupted transactions during this time;
2. The end of financial year is Esfand;
3. The Transactions of their stock, not stop more than three months;
4. The information of financial statements is available.

Considering the mentioned qualifications and using "Rahavard Novin" software, 81 companies selected for this study in this period of time.

Table 1 includes the companies' names, the price of each share and the earning per share (EPS) of each one for statistical population of this paper:

<Table 1> The data of Tehran's stock market

| Companies' name | EPS | Price | Companies' name | EPS | Price |
|---|---|---|---|---|---|
| offset | 1063 | 5697.68 | Fibr Iran | 184 | 5532.88 |
| Alborz Daru | 1195 | 5461.15 | Sar.Pardis | 405 | 1628.1 |
| Iran Tire | 1319 | 5711.27 | Sar.Tosee Sanati | 355 | 2133.55 |
| Iran transfo | 425 | 1810.5 | Sarma Afarin | 486 | 2551.5 |
| Iran Khodro | -414 | 1217.16 | Siman Orumieh | 809 | 3713.31 |
| Iran Khodro Diesel | -719 | 884.37 | Siman Tehran | 538 | 2603.92 |
| Iran Daru | 394 | 2064.56 | Siman Darud | 1263 | 4736.25 |
| Absal | 423 | 4179.24 | Siman Sepahan | 157 | 675.1 |
| Ahangari Trucktor | -446 | 3385.14 | Siman Sufian | 427 | 1648.22 |
| Behnoush | 476 | 2898.84 | Sina Daru | 2444 | 21604.96 |
| Pars Khodro | -492 | 369 | Shahd Iran | 163 | 2327.64 |
| Pars Daru | 3344 | 21033.76 | Shishe Ghazvin | -3155 | 4637.85 |
| Paksan | 671 | 2717.55 | Shishe o Gaz | 59 | 6315.36 |
| Petroshimi Abadan | 1492 | 8265.68 | Shimiaiee Sina | 498 | 4755.9 |
| Pomp Iran | 630 | 2293.2 | Faravarde Tazrighi | 1494 | 6991.92 |
| Trucktor Sazi | 446 | 2698.3 | Fanarsazi Khavar | 77 | 916.3 |
| Charkheshgar | 20 | 1115.6 | Fanarsazi Zar | -922 | 608.52 |
| Khak Chini | 2746 | 15597.28 | Fibr Iran | 184 | 5532.88 |
| Khadamat anformatik | 3332 | 14294.28 | Ghataat Otomobil | 549 | 1059.57 |
| Khorak Dam Pars | 1137 | 11426.85 | Ghand Naghshe Jahan | 740 | 9864.2 |
| Daru Abureihan | 1357 | 8304.84 | Irka Part Sanat | 249 | 1899.87 |
| Daru Osve | 2310 | 11295.9 | Kart Iran | 521 | 2360.13 |
| Daru Exir | 1283 | 6081.42 | Darupakhsh | 1333 | 5478.63 |
| Daru Amin | 696 | 2895.36 | Kashi Isfahan | -764 | 74696.28 |
| Daru Damelran | 2772 | 12501.72 | Kaghazsazi Kaveh | 788 | 4507.36 |
| Daru Zahravi | 2699 | 10256.2 | Kalsimin | 1288 | 6105.12 |
| Daru Sobhan | 834 | 3294.3 | Komak Fanar Indamin | 36 | 1268.28 |
| Daru Abidi | 404 | 2444.2 | Sar.Pardis | 405 | 1628.1 |
| Daru Farabi | 2188 | 12756.04 | Sar.Tosee Sanati | 355 | 2133.55 |
| Daru Loghman | 212 | 1358.92 | Sarma Afarin | 486 | 2551.5 |
| Daru Jaber ibn Hayyan | 659 | 2741.44 | Siman Orumieh | 809 | 3713.31 |
| | | | Siman Tehran | 538 | 2603.92 |
| Darusazi Kosar | 121 | 1395.13 | Siman Darud | 1263 | 4736.25 |
| Derakhsahn Tehran | 505 | 3873.35 | Siman Sepahan | 157 | 675.1 |
| Dashte Morghab | 152 | 2115.84 | Siman Sufian | 427 | 1648.22 |
| Radiator Iran | 30 | 222.6 | Sina Daru | 2444 | 21604.96 |
| Rikhtegari Trucktor | 130 | 850.2 | Shahd Iran | 163 | 2327.64 |
| Zamyad | -221 | 6826.69 | Shishe Ghazvin | -3155 | 4637.85 |
| Saipa | -693 | 7262.64 | Shishe o Gaz | 59 | 6315.36 |
| Saipa Azin | -21 | 2.94 | Shimiaiee Sina | 498 | 4755.9 |
| Saipa Diesel | -260 | 949 | Faravarde Tazrighi | 1494 | 6991.92 |
| Sar.Alborz | 509 | 3481.56 | Fanarsazi Khavar | 77 | 916.3 |
| Sar.Petroshimi | 155 | 4465.55 | Fanarsazi Zar | -922 | 608.52 |
| Ghataat Otomobil | 549 | 1059.57 | Ghand Naghshe Jahan | 740 | 9864.2 |
| Irka Part Sanat | 249 | 1899.87 | Kashi Isfahan | -764 | 74696.28 |
| Kart Iran | 521 | 2360.13 | Kaghazsazi Kaveh | 788 | 4507.36 |
| Darupakhsh | 1333 | 5478.63 | Kalsimin | 1288 | 6105.12 |
| Kontor sazi Iran | 2062 | 12701.92 | Komak Fanar Indamin | 36 | 1268.28 |
| Kimidaru | 227 | 2649.09 | Goruh Bahman | 425 | 1406.75 |

In following, the Knapsack mathematical model that is used for selecting the portfolio will be explained. The Knapsack model contains an objective function and one constraint. The purpose of this model is to maximize the profit of the selected portfolio subject to budget constraint.
The constraint of this model searches portfolio selection with respect to maximizing the return for investor and considering the amount of investment that he wants to invest in the stock market.

The parameters and variables of the Knapsack model are following:

$V_i$: The earning of each share of company *i* (profit of each item in the knapsack)

$W_i$: The price of each share of company *i* (price of each item in the knapsack)

$X_i$: The binary decision variable for selecting the share (0,1)

*N:* The number of studied companies (81 companies)

Considering these parameters, the mathematical model is as (1):

$$\text{Max} \sum_{i=1}^{i=n} V_i X_i$$
s.t
$$\sum_{i=1}^{i=n} W_i X_i \leq 100000$$

(1)

Model (1) insists on selecting the portfolio with the maximum profit. Also, the whole investment is 100000 units according to the assumption. This means the portfolio that will be selected cannot take more than 100000 units of investment.

International Journal on Computational Sciences & Applications (IJCSA) Vol.5, No.3, June 2015

## 4. IMPLEMENTING THE COA MODEL TO PROPOSED MODEL

In this part, we want to implement the COA (version of [15]) to Knapsack model of portfolio selection problem. The proposed mathematical model coded in MATLAB. The constraint of the model related to the amount of investment, applied by using the multiplicative violation function and the violation factor $a=10$ in the objective function.

According to the Violation function, the constraint will add to objective function as (2).

$$Violation = Max\left(\frac{\sum_{i=1}^{i=n} WiXi}{100000} - 1, 0\right) \quad (2)$$

As can be seen the value of violation function is zero or a positive magnitude. If the constraint of model is satisfied, the violation magnitude becomes zero and the violation function doesn't affect the objective function. The objective function with the violation is as (3).

$$Max \sum_{i=1}^{i=n} ViXi \, (1+\alpha(Violation)) \quad (3)$$

One of the most important issues in meta-heuristic algorithms is parameter setting. The COA parameters are as following:

Number of initial population=10, Minimum number of eggs for each cuckoo=2, Maximum number of eggs for each cuckoo=4, Maximum iterations of the COA=200, Maximum number of cuckoos that can live at the same time=10.

Portfolio selection problem was solved by the COA in several times. Some of obtained solutions are shown in Table 2. Also, the convergence rate of COA is shown in Figure 3.

<Table2> The results of COA

| Iteration number | Elapsed time | Objective function |
|---|---|---|
| 1 | 8.92 | 35672.00 |
| 2 | 8.8 | 36041.00 |
| 3 | 9.1 | 35812.00 |

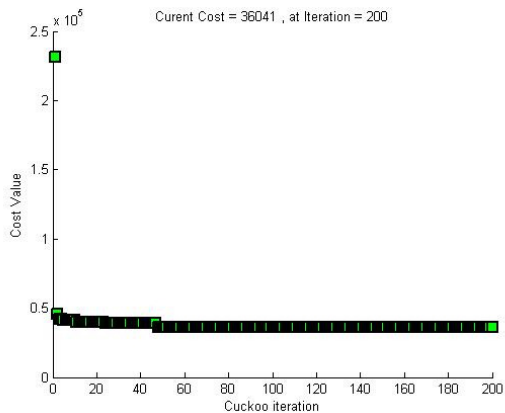

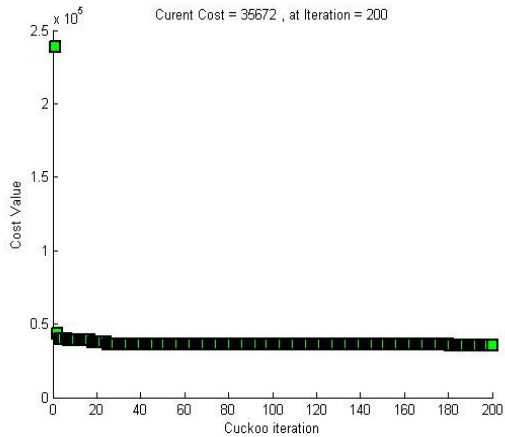

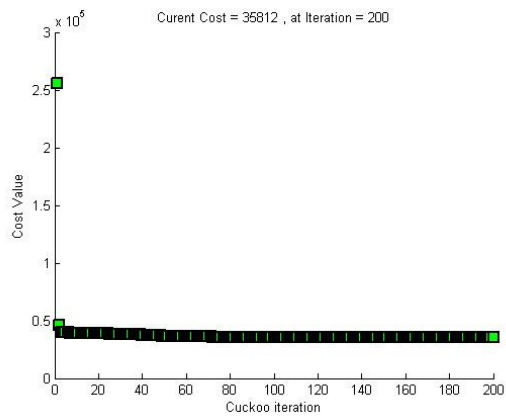

<Figure 3> The convergence rate of COA

Also, portfolio selection problem is optimized by the Genetic algorithm using the parameters that listed in the following and the results are shown in Table 3. The convergence rate of the Genetic algorithm is shown in Figure 4. Also, the version of this Genetic algorithm is GA of MATLAB R2014b software.

The parameters of Genetic algorithm include Population size=150, Crossover percentage=0.9, Mutation percentage=0.3, Mutation rate=0.02.

<Table 3> The results of Genetic Algorithm

| Number of iteration | Elapsed time | Objective function |
|---|---|---|
| 1 | 10.1 | 32691.00 |
| 2 | 9.8 | 32522.00 |
| 3 | 9.6 | 32503.00 |

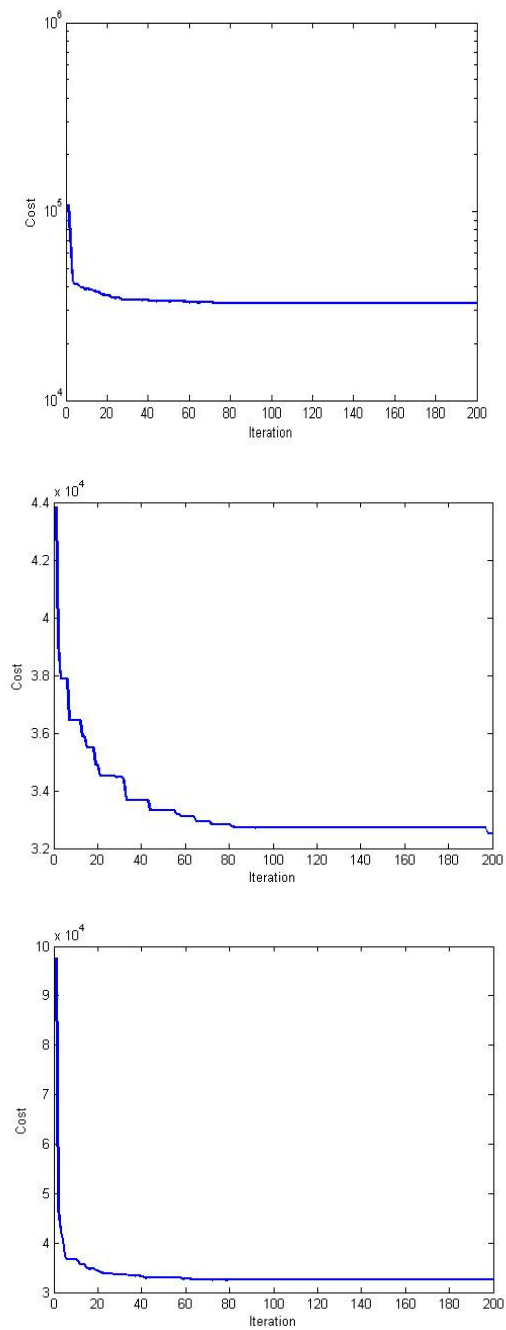

<Figure 4> The convergence rate of Genetic Algorithm

As Table 2 and 3 show that the COA obtains better results in the less elapsed time than Genetic algorithm. Also, considering convergence rate, the COA obtains solution in less iteration than the Genetic algorithm and this issue represents that the speed of the COA higher than the Genetic Algorithm.

## 5.CONCLUSION

In this paper, we tried to select the best portfolio with respect to maximizing the profit using the data of Tehran's stock market.

For achieving this purpose, first the mathematical model (the Knapsack model) of portfolio selection problem that obtains the most possible profit was modelled.

Assuming the data of Tehran's stock market and the Knapsack model, the COA applied to select the optimized portfolio of stocks. Then, this problem implemented in the Genetic Algorithm too and the portfolio selected this way. At last, two portfolios obtained (from the COA and the Genetic Algorithm). Comparing the final results showed more convergence rate and accuracy of the COA rather than the Genetic Algorithm in low iteration.

For continuing the research, it is possible to expand the model and development an integer model considering the risk factor.